\documentclass[12pt]{amsart}
\usepackage{hyperref}
\hypersetup{nesting=true,debug=true,naturalnames=true}
\usepackage{graphicx,amssymb}

\usepackage[square,sort&compress,comma,numbers]{natbib}
\usepackage{nicefrac,xcolor,upref,version}

\usepackage{amscd,amsthm,amsmath,amssymb,amsfonts}

\usepackage{mathrsfs}

\newtheorem{theorem}{Theorem}[section]
 
\newtheorem{proposition}[theorem]{Proposition}

\numberwithin{equation}{section}

\def\C{{\mathbb C}} 
\def\eps{{\varepsilon}}

\def\tig{{\widetilde g}}
\def\cL{{\mathcal L}}

\def\tif{{g}}

\def\cFrac#1#2{%
  \begin{array}{@{}c@{}}\multicolumn{1}{c|}{#1}\\%
  \hline\multicolumn{1}{|c}{#2}\end{array}\;}

\allowdisplaybreaks

\begin{document}

\title[Approximation to Thue--Morse rational numbers]{On the rational approximation 
to Thue--Morse rational numbers}

\author{Yann Bugeaud}
\address{Universit\'e de Strasbourg, Math\'ematiques,
7, rue Ren\'e Descartes, 67084 Strasbourg  (France)}
\email{bugeaud@math.unistra.fr}

\author{Guo-Niu Han}
\address{Universit\'e de Strasbourg, Math\'ematiques,
7, rue Ren\'e Descartes, 67084 Strasbourg  (France)}
\email{guoniu@gmail.com}

\begin{abstract}
Let $b \ge 2$ and $\ell \ge 1$ be integers. 
We establish that there is an absolute real number $K$ such that all the partial quotients of 
the rational number 
$$
\prod_{h = 0}^\ell \, (1 - b^{-2^h}),
$$
of denominator $b^{2^{\ell+1} - 1}$, 
do not exceed $\exp(K (\log b)^2 \sqrt{\ell} 2^{\ell/2})$. 
\end{abstract}


\maketitle

\section{Introduction}\label{intro}

An easy covering argument which goes back to Cantelli shows that, for almost all real numbers 
$\xi$ (with respect
to the Lebesgue measure) and for every positive $\eps$, the inequality
$$
\biggl| \xi - {p \over q} \biggr| > {1 \over q^{2 + \eps}}
$$
holds for every sufficiently large $q$. 
However, it is often a very difficult problem to show that a given real number shares 
this property, unless its continued fraction expansion is explicitly determined. 
This is known to be the case for any irrational real algebraic number, by Roth's theorem, and 
for only a few other real numbers defined by
their expansion in some integer base. 
Let 
$
{\bf t} = t_0 t_1 t_2 \ldots
$
denote the Thue--Morse word over $\{-1, 1\}$ defined by $t_0 = 1$,
$t_{2k} = t_k$ and $t_{2k + 1} = - t_k$ for $k \ge 0$. Then, the Thue--Morse 
generating series
$\xi_{{\bf t}} (z)$ is given by 
\begin{align*}
\xi_{{\bf t}} (z) & = \sum_{k \ge 0} \, t_{k} z^{-k}  
=   1 - z^{-1} - z^{-2} + z^{-3} - z^{-4} + z^{-5} + z^{-6} -  \ldots   \\
& = \prod_{h \ge 0} \, (1 - z^{-2^h}). 
\end{align*}
By means of a non-vanishing result obtained in \cite{APWW98} for the Hankel determinants
associated with the Thue--Morse sequence,   
Bugeaud \cite{Bu11} established that, for any given positive $\eps$ and any integer $b \ge 2$, the 
Thue--Morse--Mahler number 
$$
\xi_{{\bf t}} (b) =  \sum_{k \ge 0} \, {t_{k} \over b^{k}}
=  1 - {1 \over b} - {1 \over b^2} + {1 \over b^3} - {1 \over b^4} + {1 \over b^5} + {1 \over b^6} - {1 \over b^7} -
{1 \over b^{8}} + \ldots 
$$
satisfies the inequality 
$$
\Bigl| \xi_{{\bf t}} (b) - {p \over q} \Bigr| > {1 \over q^{2 + \eps}}, 
$$
for every sufficiently large $q$. 
Subsequently, his result has been considerably improved by 
Badziahin and Zorin \cite{BaZo20}, who showed that there 
exists a positive real number $K$ such that the stronger inequality 
\begin{equation} \label{eq:1.1}
\Bigl| \xi_{{\bf t}} (b) - {p \over q} \Bigr| > {1 \over q^2 \exp( K \log b \, \sqrt {\log q \, \log \log q})}  
\end{equation} 
holds as soon as $q$ is large enough. Thus, all the partial quotients of 
$\xi_{{\bf t}} (b)$ are rather small. Note that, in view of \cite[Th. 11]{BaZo20}, the number $K$ occurring 
in \cite[Th. 2]{BaZo20} must depend on $b$ and can be taken equal to
$\log b$ times some number depending only 
on the series $f(z)$ occurring in the statement of \cite[Th. 2]{BaZo20}.

Observe that the Thue--Morse power series $\xi_{{\bf t}} (z)$
is the limit of the sequence of rational functions
$$
f_\ell (z) = \prod_{h = 0}^\ell \, (1 - z^{-2^h}). 
$$
More precisely, we have
$$
\xi_{{\bf t}} (z)  = f_\ell (z) + O (z^{-2^{\ell +1}}), \quad \ell \ge 1,   
$$
and
\begin{equation} \label{eq:1.2}
|\xi_{{\bf t}} (x) - f_\ell (x) | \le \frac{1}{(|x| - 1) |x|^{2^{\ell + 1} - 1}}, \quad \ell \ge 1, 
x \in \C, |x| > 1.   
\end{equation} 
Let $b \ge 2$ and $\ell \ge 1$ be integers. For a rational number $p/q$, we derive from \eqref{eq:1.2} that 
\begin{equation} \label{eq:1.3}
\biggl| \, \Bigl| \xi_{{\bf t}} (b) - \frac{p}{q} \Bigr| - \Bigl| f_\ell (b) - \frac{p}{q} \Bigr| \, \biggr|
\le \frac{1}{(b - 1) b^{2^{\ell +1} - 1}}.   
\end{equation} 
Consequently, $\xi_{{\bf t}} (b)$ and $f_\ell (b)$ have the same 
first partial quotients. 
To see this, let $p_n/q_n$ be the convergent to $\xi_{{\bf t}} (b)$ with 
$q_n \le b^{2^\ell}$ and $n$ maximal for this property. 
We assume that $\ell$ is sufficiently large to ensure that $n \ge 8$. 
A short calculation shows that 
\begin{equation} \label{eq:1.4}
q_n \ge q_{n-1} + q_{n-2} \ge \ldots \ge 8 q_{n-5}.   
\end{equation}  
By a result of Borel \cite[Ch. I, Th. 5B]{SchmLN}, there exists $\eps$ in $\{0, 1, 2\}$ such that 
$$
\Bigl| \xi_{{\bf t}} (b) - \frac{p_{n-5-\eps}}{q_{n-5-\eps}} \Bigr| \le \frac{1}{\sqrt{5} q^2_{n-5- \eps}}. 
$$
It then follows from \eqref{eq:1.3} and \eqref{eq:1.4} that 
$$
\Bigl| f_\ell (b) - \frac{p_{n-5-\eps}}{q_{n-5-\eps}} \Bigr| \le \frac{1}{\sqrt{5} q^2_{n-5- \eps}} 
+ \frac{2}{q_n^2} \le \Bigl( \frac{1}{\sqrt{5}} + \frac{1}{32} \Bigr)  \frac{1}{q^2_{n-5- \eps}} 
< \frac{1}{2 q^2_{n-5- \eps}},
$$
which, by a classical theorem of Legendre \cite[Ch. I, Th. 5C]{SchmLN}, 
implies that $p_{n-5-\eps} / q_{n-5-\eps}$ is a convergent 
of $f_\ell (b)$. Consequently, $\xi_{{\bf t}} (b)$ 
and $f_\ell (b)$ have the same $n-7$ first partial quotients. 
By \eqref{eq:1.1}, these
partial quotients are rather small. 
However, \eqref{eq:1.1} gives no information on the remaining partial quotients 
of $f_\ell (b)$, thus, in particular, on the rate with which the rational number $f_\ell (b)$ 
of denominator $b^{2^{\ell +1} - 1}$ is approximated by rational numbers $p/q$ 
of denominator $q$ greater than $b^{2^{\ell }}$. 
In the present note, we address this question and show that an inequality 
like \eqref{eq:1.1} remains true for every convergent of $f_\ell (b)$.

\begin{theorem} \label{th:main} 
There exists a positive real number $K$ such that, 
for every integer $b \ge 2$ and every integer $\ell \ge 2$, the inequality
$$
\Bigl| \prod_{h = 0}^\ell \, (1 - b^{-2^h}) - {p \over q} \Bigr| 
> {1 \over q^2 \exp( K \log b \, \sqrt {\log q \, \log \log q})},
$$
holds for every rational number $p/q$ different from $f_\ell (b)$. Write
$$
f_\ell (b) = \prod_{h = 0}^\ell \, (1 - b^{-2^h})  =  [0; a_1^{(\ell)}, a_2^{(\ell)}, \ldots , a_{\cL(\ell)}^{(\ell)}],
$$
with $a_{\cL(\ell)}^{(\ell)} \ge 2$. 
The partial quotients $a_j^{(\ell)}$ of $f_\ell (b)$ are all at most equal to 
$b^{2 K \sqrt{\ell} \sqrt{ \log b \, \log \log 3b} 2^{\ell / 2}}$. 
There exists a positive real number $C$, depending only on $b$, such that the 
length $\cL(\ell)$ of the continued fraction of $f_\ell (b)$ exceeds $C 2^{\ell/2} / \sqrt{\ell}$. 
\end{theorem}

The second assertion of Theorem \ref{th:main} immediately follows from the 
first one and the theory of continued fractions. The last assertion already follows from \eqref{eq:1.1}. 

Theorem \ref{th:main} is mainly motivated by the very few known results 
on continued fraction expansions of sequences of rational numbers. 
Pourchet \cite{Po72} (see also \cite{vdP81,Bu09}) 
proved that, for all coprime integers $a$ and $b$ with $1 < b < a$ and for every positive $\eps$,
there exists a positive $C$, depending only on $\eps$ and on the prime divisors of $a$ and $b$, such that 
all the partial quotients of $(a/b)^n$ are less than $C b^{\eps n}$.
This was subsequently extended to quotients of power sums by Corvaja and Zannier \cite{CZ05}, 
with a similar conclusion. Consequently, the length of the continued fraction 
expansion of $(a/b)^n$ (resp., of $(a^n - 1) / (b^n - 1)$) tends to infinity with $n$. 
We stress that the conclusion of Theorem \ref{th:main} is much stronger. 

\begin{figure}
  \centering
  \includegraphics[width=0.5\textwidth]{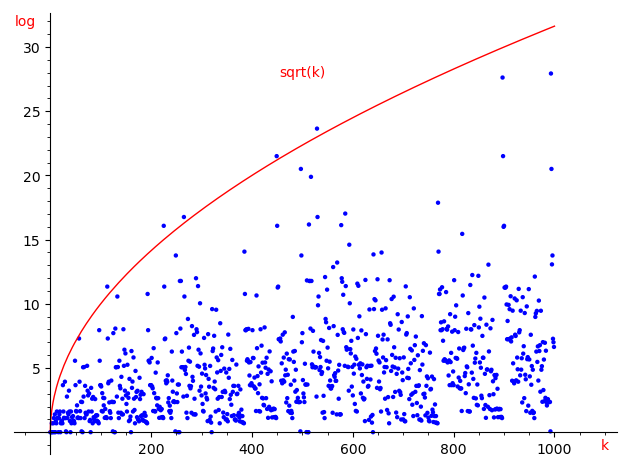}\includegraphics[width=0.5\textwidth]{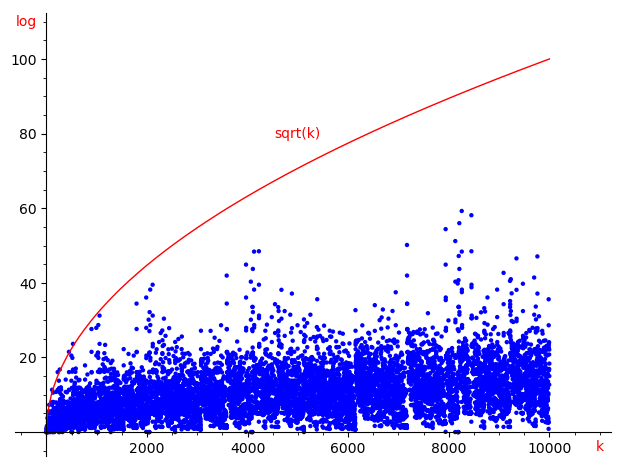}
  \caption{Logarithm of the absolute values of the coefficients of the $k$-th convergent}
  \label{figure_F12}
\end{figure}

The function $q \mapsto \exp ( \sqrt{ \log q \log \log q})$ occurring in \eqref{eq:1.1} is a consequence 
of the bound of order $(c_1 k)^{c_2 k}$ obtained in \cite{BaZo20} 
for the absolute values of the coefficients of the numerator 
and denominator of the $k$-th convergent to $\xi_{{\bf t}} (z)$. 
Numerical experiments suggest that a better bound of the shape $c_3^{\sqrt{k}}$ should hold
(here, $c_1, c_2$ and $c_3$ are absolute, positive real numbers). Such a result 
seems difficult to establish. See Figure \ref{figure_F12}.

To prove \eqref{eq:1.1}, Badziahin and Zorin \cite{BaZo20} used 
that all the partial quotients of the continued fraction expansion 
of $\xi_{{\bf t}} (z)$ are linear, a result established by Badziahin \cite{Ba19}. 
Here, we first show that  all the partial quotients of the rational functions 
$f_\ell (z)$, $\ell \ge 1$, are linear. This is the main novelty 
of the present note and the object of Section \ref{sec:2}. 
Then, in Section \ref{sec:3}, we establish Theorem \ref{th:main} by adapting 
to our purpose the argumentation of \cite{BaZo20}. Finally, in the last section, we 
discuss another example.

\section{The partial quotients of the rational functions $f_\ell (z)$} \label{sec:2}

For a non-zero rational number $x$, we let $\nu_2 (x)$ denote 
its $2$-adic valuation, that is, the exponent of $2$ 
in its decomposition in product of prime factors. Put $\nu_2 (0) = + \infty$.

\begin{proposition}\label{Prop:nonzero}
Let $v^{(\ell)}_j$ $(\ell\geq 0, 0\leq j\leq 2^{\ell+1}-1)$ be the family of rational constants defined by
$$
v^{(0)}_0=1,\quad
v^{(0)}_1=1,
$$
and, for $\ell\geq 1$,
\begin{align*}
	v^{(\ell)}_0&=1, \\
	v^{(\ell)}_1 &=2, \\
	v^{(\ell)}_{2j} &= -\frac{ v^{(\ell-1)}_{j}}{    v^{(\ell)}_{2j-1} } ,\quad 1\leq j\leq 2^{\ell}-1, \\
	v^{(\ell)}_{2j+1} &= 1+(-1)^j   -   v^{(\ell)}_{2j} , \quad 1\leq j\leq 2^{\ell}-1.
\end{align*}
Then,	$\nu_2 (v^{(\ell)}_j)$  equals to

$\bullet$ $1$ for $\ell\geq1$ and $j=1$;

$\bullet$ $-1$ for $\ell\geq 1$ and $j=2^\ell$ or $j=2^\ell+1$;

$\bullet$ $0$ otherwise.

\noindent As a consequence, we have $v^{(\ell)}_j\not=0$ for all 
$\ell\geq 0, 0\leq j\leq 2^{\ell+1}-1$.
\end{proposition}

The first values of $v^{(\ell)}_j$ are given in the following table:
\begin{equation*}
  \begin{array}{*{19}c}
    \ell \setminus j & 0 & 1 & 2& 3& 4& 5& 6& 7& 8& 9& 10 & 11& 12& 13& 14& 15& \\
0&1 & 1 \\
1&1 & 2 &-\frac{1}{2} & \frac{1}{2} & \\
2&1 & 2 & -1 & 1 &\frac{1}{2} & \frac{3}{2} & -\frac{1}{3} & \frac{1}{3} \\
3&1& 2& -1& 1& 1& 1& -1& 1& -\frac{1}{2}& \frac{5}{2}& -\frac{3}{5}& \frac{3}{5}& \frac{5}{9}& \frac{13}{9}& -\frac{3}{13}& \frac{3}{13}\
\end{array}%
\end{equation*}

Consequently, the first values of $\nu_2 (v^{(\ell)}_j)$ are given in the following table:
\begin{equation*}
  \begin{array}{*{19}c}
    \ell \setminus j & 0 & 1 & 2& 3& 4& 5& 6& 7& 8& 9& 10 & 11& 12& 13& 14& 15& \\
0& 0 & 0 \\
1&0 & 1 & -1 & -1 & \\
2&0 & 1 & 0 & 0 & -1  &-1 & 0 & 0 \\
3&0& 1 & 0& 0& 0& 0& 0& 0& - 1 &-1 & 0& 0& 0& 0& 0& 0\\
\end{array}%
\end{equation*}

\begin{proof}
We proceed by induction on $\ell$. 
Recall that, for any non-zero rational numbers $x, y$, we have $\nu_2 (x/y) = \nu_2 (x) - \nu_2 (y)$ 
and $\nu_2 (x + y) \ge \min\{\nu_2(x), \nu_2(y)\}$, with equality if $\nu_2(x) \not= \nu_2(y)$. 
The tables above show that the proposition holds for $0  \le \ell \le 3$. 
Let $\ell \ge 4$ be an integer such that the proposition holds for $\ell-1$.
By definition, we have
$$
\nu_2 (v^{(\ell)}_0) = 0, \quad \nu_2 (v^{(\ell)}_1) = 1, \quad
\nu_2 (v^{(\ell)}_2) = \nu_2 (-1) =0.
$$
Since $\nu_2 (v^{(\ell)}_{2j+1}) = \nu_2 (v^{(\ell)}_{2j})$ for $j = 1, \ldots , 2^\ell - 1$ such that 
$\nu_2 (v^{(\ell)}_{2j}) \not= 1$, we derive that $\nu_2 (v^{(\ell)}_3) = 0$, thus, $v^{(\ell)}_3$ is nonzero 
and $\nu_2 (v^{(\ell)}_4) = 0$. Continuing in this way, we get that 
$v^{(\ell)}_5, \ldots , v^{(\ell)}_{2^\ell - 1}$ are all nonzero and
$$
 \nu_2 (v^{(\ell)}_5) = \ldots = \nu_2 (v^{(\ell)}_{2^\ell - 1}) = 0, \quad 
\nu_2 (v^{(\ell)}_{2^\ell }) = \nu_2 (v^{(\ell-1)}_{2^{\ell-1 }}) = -1.
$$
Then, $\nu_2 (v^{(\ell)}_{2^\ell +1}) = - 1$. Thus, $v^{(\ell)}_{2^\ell +1}$ is nonzero 
and $\nu_2 (v^{(\ell)}_{2^\ell + 2}) = -1 - (-1) = 0$. We derive that $\nu_2 (v^{(\ell)}_{2^\ell + 3}) = 0$, 
thus, $v^{(\ell)}_{2^\ell + 3}$ is nonzero 
and $\nu_2 (v^{(\ell)}_{2^\ell + 4}) = 0$. Continuing in this way, we get that 
$v^{(\ell)}_{2^\ell + 5}, \ldots , v^{(\ell)}_{2^{\ell +1} - 1}$ are all nonzero and
$$
\nu_2 (v^{(\ell)}_{2^\ell + 5}) = \ldots = \nu_2 (v^{(\ell)}_{2^{\ell+1} - 1}) = 0.
$$
This completes the induction step. 
\end{proof}

Set 
$$
\tif_\ell(z)=\frac{1}{z} f_\ell(z) = \frac{1}{z} \prod_{h = 0}^\ell \, (1 - z^{-2^h}) 
$$
and
$$
\tif_\ell(z)= [0; a_1(z), a_2(z), \ldots, a_m(z)]
=\cFrac{1}{a_1(z)}+\cFrac{1}{a_2(z)} + \cdots + \cFrac{1}{a_m(z)}.
$$
where $a_i(z)$ is in $\mathbb{Q}[z]$ for $1\leq i\leq m$.
The following theorem, which can be seen as a finite version of \cite[Prop. 3.3]{BaZo15}, 
shows that all the $a_i(z)$ are polynomials of degree one.

\begin{theorem}\label{th:fzfrac}
	Let $v^{(\ell)}_j$ $(\ell\geq 0, 0\leq j\leq 2^{\ell+1}-1)$ be the family of rational numbers 
	defined in Proposition \ref{Prop:nonzero}. Then,
\begin{equation}\label{eq:fzfrac}
\tif_\ell(z)=
	\cFrac{v^{(\ell)}_0  }{ z +  1} 
	+ \cFrac{v^{(\ell)}_1  }{ z -1} 
	+ \cFrac{v^{(\ell)}_2  }{ z +1} 
	+ \cdots
	+ \cFrac{v^{(\ell)}_{2^{\ell+1}-1}  }{ z - 1}, \quad \ell \ge 0. 
\end{equation}
For $\ell \ge 0$, all the partial quotients in the continued fraction expansion of $\tif_\ell(z)$ are 
polynomials of degree one. 
\end{theorem}

The last assertion of Theorem \ref{th:fzfrac} immediately follows from 
Proposition~\ref{Prop:nonzero}.  

\begin{proof}
	We prove identity \eqref{eq:fzfrac} by induction on $\ell$.
Since
$$
\frac{1}{z}(1 - z^{-1}) = \cFrac{1}{z + 1} 
+ \cFrac{1}{z-1 } 
$$
and
$$
\frac{1}{z}(1 - z^{-1})(1 - z^{-2}) = \cFrac{1}{z + 1} 
+ \cFrac{2}{z-1 } - \cFrac{1/2}{ z+1} + \cFrac{1/2}{ z -1},
$$
identity \eqref{eq:fzfrac} is true for $\ell=0,1$. 
	Let $k \ge 1$ be an integer and suppose that \eqref{eq:fzfrac} is true for $\ell\leq k$.
Set
\begin{equation}\label{def:gz}
h_{k+1}(z)=
	\cFrac{v^{({k+1})}_0  }{ z +  1} 
	+ \cFrac{v^{({k+1})}_1  }{ z -1} 
	+ \cFrac{v^{({k+1})}_2  }{ z +1} 
	+ \cdots
	+ \cFrac{v^{({k+1})}_{2^{{k+2}}-1}  }{ z - 1}.
\end{equation}
It suffices prove that $\tif_{k+1}(z)=h_{k+1}(z)$.
	From the even contraction theorem (see, for example, \cite[Theorem 2.1(1)]{Han20}), 
	we have
\begin{equation}\label{eq:GCF}
	h_{k+1}(z)=
b_0 +  \cFrac{a_1}{b_1} 
+ \cFrac{a_2}{b_2} 
+ \cFrac{a_3}{ b_3} 
+ {\cdots}, 
\end{equation}
where
\begin{align*}
	a_1 &=  v^{(k+1)}_0 (z-1),   \cr
	a_2 &=  - v^{(k+1)}_1 v^{(k+1)}_2 { (z-1)}, \cr
	a_j &=  -  v^{(k+1)}_{2j-3} v^{(k+1)}_{2j-2}  (z-1)^2, \quad 3\leq j\leq 2^{k+1}, \cr
	a_j &=  0, \quad j>2^{k+1}, \cr
	b_0  &=  0, \cr
	b_1 &=  (z+1)(z-1) +v^{(k+1)}_1, \\
	b_j &= (z-1)  \bigl(  (z-1)(z+1)  +  v^{(k+1)}_{2j-1}  +   v^{(k+1)}_{2j-2}  \bigr), \quad 2\leq j\leq 2^{k+1}, \\
	b_j &=  0, \quad j>2^{k+1}. 
\end{align*}
By removing the common factors in numerators and denominators, we obtain
\begin{equation}
	h_{k+1}(z)=
b_0 +  \cFrac{a_1}{b_1} 
+ \cFrac{a_2}{b_2} 
+ \cFrac{a_3}{ b_3} 
+ {\cdots}
\end{equation}
where
\begin{align*}
	a_1 &=  v^{(k+1)}_0 (z-1),   \cr
	a_2 &=  - v^{(k+1)}_1 v^{(k+1)}_2 , \cr
	a_j &=  -  v^{(k+1)}_{2j-3} v^{(k+1)}_{2j-2}, \quad 3\leq j\leq 2^{k+1}, \cr
	a_j &=  0, \quad j>2^{k+1}, \cr
	b_0  &=  0, \cr
	b_1 &=  (z+1)(z-1) +v^{(k+1)}_1, \\
	b_j &=     (z-1)(z+1)  +  v^{(k+1)}_{2j-1}  +   v^{(k+1)}_{2j-2} , \quad 2\leq j\leq 2^{k+1}, \\
	b_j &=  0, \quad j>2^{k+1}. 
\end{align*}
Using the recurrence relations defined in the statement of Theorem \ref{th:fzfrac},
a quick calculation shows that we have
\begin{equation*}
h_{k+1}(z)=
	\cFrac{z-1 }{ z^2 +  1} 
	+ \cFrac{v^{({k})}_1  }{ z^2 -1} 
	+ \cFrac{v^{({k})}_2  }{ z^2 +1} 
	+ \cdots
	+ \cFrac{v^{({k})}_{2^{{k+1}}-1}  }{ z^2 - 1}.
\end{equation*}
This implies that 
$ h_{k+1}(z) = (z-1) \tif_{k}(z^2) = \tif_{k+1}(z) $.
\end{proof}
For a more general $(z+1, z-1)$ phenomenon, see \cite[Lemma 3.1]{Han16}.

\section{Completion of the proof of Theorem \ref{th:main}}  \label{sec:3}

The key new ingredient for the proof of Theorem \ref{th:main} is the fact that 
all the partial quotients of the rational functions $f_\ell (z)$ are linear. This allows us to follow the 
argumentation of \cite{BaZo20}, with some minor changes. 
For the sake of readability, we keep most of the notation of \cite{BaZo20} and we 
sketch how to adapt the proof of \cite[Th. 11]{BaZo20}. 
Instead of working with the (infinite) power series $g_{{\bf u}} (z)$, we fix a positive integer $\ell$
and work with the (finite) power series 
$$
g_\ell (z) = z^{-1} \, f_\ell (z) = z^{-1} - z^{-2} - z^{-3} + z^{-4} + \ldots + (-1)^{\ell} z^{-2^\ell}. 
$$
Since, by Theorem \ref{th:fzfrac}, all 
the partial quotients of $g_\ell (z)$ are linear, the auxiliary results in \cite{BaZo20} hold. 
Furthermore, for $m \ge 1$, we have 
\begin{align*}
g_{\ell + m}  (z) & = z^{-1} \, \prod_{h = 0}^{\ell + m} \, (1 - z^{-2^{h}}) \\
& = g_\ell (z^{2^m}) \, z^{{2^m} - 1} \, \prod_{h = 0}^{m-1} \, (1 - z^{-2^{h }}) 
 = g_\ell (z^{2^m}) \, \prod_{h = 0}^{m-1} \, (  z^{2^{h }} - 1). 
\end{align*}
Denoting by $p_{k, \ell} (z) / q_{k, \ell} (z)$ the convergents to $q_\ell (z)$, where $k=1, \ldots , 2^{\ell +1}$, 
and defining $p_{k, \ell, m} (z)$ and $q_{k, \ell, m} (z)$ as in \cite[(2.20)]{BaZo20}, the  
analogue of \cite[(2.25)]{BaZo20} holds, namely, we have 
$$
\Bigl| g_{\ell + m}  (b) - {p_{k,\ell, m} (b) \over q_{k,\ell, m} (b)} \Bigr| 
\le {2 (k+1) k^{k/2} 2^m \over q_{k, \ell, m} b^{k 2^m + 1}}, 
$$
for $k = 1, \ldots , 2^{\ell +1}$. By \cite[Lemma 9]{BaZo20}, the integer $q_{k, \ell, m}$ is controlled 
and is comparable to $b^{k 2^m}$. 

Take now a large integer $L$ (which corresponds to the integer $\ell$ in the statement of the theorem)
and study the rate of approximation to the 
rational number $g_L (b)$ of denominator $b^{2^{L+1}}$ by rational numbers $p/q$ 
of denominator less than $b^{2^{L+1}}$. 
We follow the argument of \cite{BaZo20} and look for a power of $b$ 
close to $q$. We use the fact that every integer $n$ less than $2^{L+1}$ is rather close 
to a product $k 2^m$, where $L = m + \ell$ and $k \le 2^{\ell +1 }$. 
The latter constraint comes from the construction of our finite sequence of good rational 
approximations to $g_L (b)$. It is not required to hold in \cite{BaZo20}.

Let $p/q$ be a rational number with $q < b^{2^{L+1}}$ and
$q$ sufficiently large (it is sufficient to assume that $q$ exceeds $b^{\kappa_1}$, 
for some absolute constant $\kappa_1$) 
to guarantee that the real number $x$
defined in \cite[(3.2)]{BaZo20} satisfies $x > b^2$.  

As in \cite[(3.8)]{BaZo20}, set $t = (\log x) / (\log b)$. Note that $t > 2$. 
Let $\tau \ge 2$ be a real number. 
Denote by $\log_2$ the logarithm in base $2$. Assume that  
\begin{equation}  \label{eq:3.1}
t + 2  \tau \sqrt {t \log_2 t} < 2^{L+1}.      
\end{equation}
This inequality holds if for a suitable positive $\kappa_2$, depending only on $\tau$, we have
\begin{equation}  \label{eq:3.2}
q < \exp \bigl( - \kappa_2 \log b \, \sqrt{\log_2 q \log_2 \log_2 q} \, \bigr) b^{2^{L+1}}.   
\end{equation}
There exist integers $n$ and $m$ and a real number $\alpha$ such that $2^m$ divides $n$ and
$$
t \le n \le t + 2  \tau \sqrt {t \log_2 t}, \quad 
2^m = \alpha \tau \sqrt {t \log_2 t}, \quad 1 \le \alpha \le 2. 
$$
It follows from \eqref{eq:3.1} that the integer 
$$
k = {n \over 2^m} \le {t \over  \alpha \tau \sqrt {t \log_2 t}} + {2 \over \alpha} 
$$
is less than $2^{L-m+1}$, as required. 

To use the results of \cite{BaZo20}, we also have 
to check \cite[(2.26)]{BaZo20} and \cite[(2.34)]{BaZo20}, that is, the inequalities
$$
2^{2^m} \ge 4 (k+1) k^{k/2}, \quad 2^{2^m} > 3 k^{k/2}. 
$$
Since 
$$
\alpha \tau \sqrt {t \log_2 t} \ge \Bigl( {\sqrt{t} \over  \alpha \tau \sqrt {\log_2 t}} + {2 \over \alpha}  \Bigr) 
\log \Bigl( {\sqrt{t} \over  \alpha \tau \sqrt {\log_2 t}} + {2 \over \alpha}  \Bigr)
$$
holds if $\tau$ is large enough, both inequalities are satisfied if $\tau$ is large enough. 

We conclude that we have an inequality similar to \cite[(3.1)]{BaZo20} provided that 
$q$ exceeds $b^{\kappa_1}$ and satisfies \eqref{eq:3.2}. This proves the theorem. 

\section{A further example}  \label{sec:4} 

The method developed in Section \ref{sec:2} is not specific to the Thue--Morse sequence and 
may be used to derive a similar conclusion for other sequences. In this section, we give 
a further example, whose corresponding infinite product was considered in \cite{Ba19,BaZo20}.

For nonzero integers $u, v$ with $u^2 \not= v$, consider
$$
\tig_{u, v, \ell} (z)=\frac{1}{z} \prod_{h=0}^\ell (1+uz^{-3^h}+vz^{-2\cdot 3^h}), \quad \ell \ge 0.
$$
Following \cite[Theorem 1.2]{Ba19}, we define
\begin{align*}
\alpha^{(0)}_1 &= -u,\quad
	\alpha^{(0)}_2 =\frac{u^3 - 2uv}{u^2 - v},\quad
	\alpha^{(0)}_3= \frac{uv}{u^2 - v},\\
\beta^{(0)}_1&=1,\quad
\beta^{(0)}_2=u^2-v,\quad
	\beta^{(0)}_3=\frac{v^3}{u^4 - 2u^2v + v^2},\\
\alpha^{(\ell)}_1&=-u,\quad
	\alpha^{(\ell)}_2=\frac{u^3 - 2uv+u}{u^2 - v},\quad
	\alpha^{(\ell)}_3= \frac{uv - u}{u^2 - v}, \\
\beta^{(\ell)}_1&=1,\quad
\beta^{(\ell)}_2=u^2-v,\quad
	\beta^{(\ell)}_3=\frac{u^4 - 3u^2v + v^3 + u^2}{u^4 - 2u^2v + v^2},\\
\alpha^{(\ell)}_{3k+4} &= -u,\quad
	\beta^{(\ell)}_{3k+4} = \frac{\beta^{(\ell-1)}_{k+2}}{\beta^{(\ell)}_{3k+3}\beta^{(\ell)}_{3k+2}},\quad
	\beta^{(\ell)}_{3k+5}=u^2-v-\beta^{(\ell)}_{3k+4},\\
	\alpha^{(\ell)}_{3k+5} &= u- \frac{\alpha^{(\ell-1)}_{k+2}+uv-\alpha^{(\ell)}_{3k+2}\beta^{(\ell)}_{3k+4}}{\beta^{(\ell)}_{3k+5}},\quad
	\alpha^{(\ell)}_{3k+6}=u-\alpha^{(\ell)}_{3k+5},\\
	\beta^{(\ell)}_{3k+6} &= v-\alpha^{(\ell)}_{3k+5} \alpha^{(\ell)}_{3k+6}.
\end{align*}
We claim that 
$$
\tig_{u, v, \ell} (z)=
	\cFrac{\beta^{(\ell)}_1  }{ z +  \alpha^{(\ell)}_1} 
	+ \cFrac{\beta^{(\ell)}_2  }{ z + \alpha^{(\ell)}_2} 
	+ \cFrac{\beta^{(\ell)}_3  }{ z + \alpha^{(\ell)}_3} 
	+ \cdots
	+ \cFrac{\beta^{(\ell)}_{3^{\ell+1}}  }{ z + \alpha^{(\ell)}_{3^{\ell+1}}}, \quad \ell \ge 0,
$$
where some of the rational numbers $\beta^{(\ell)}_{j}$ may vanish. 

In the sequel, we consider only the case $u=v=-1$, that is,
$$
\tig_\ell(z)=\frac{1}{z} \prod_{h=0}^\ell (1-z^{-3^h}-z^{-2\cdot 3^h}),
$$
and establish the non-vanishing result, also by using the $2$-adic valuation. Set
\begin{align*}
\alpha^{(0)}_1 &= 1,\quad
	\alpha^{(0)}_2 =-\frac{3}{2},\quad
	\alpha^{(0)}_3= \frac{1}{2},\\
\beta^{(0)}_1&=1,\quad
\beta^{(0)}_2=2,\quad
	\beta^{(0)}_3=-\frac{1}{4},\\
\alpha^{(\ell)}_1&=1,\quad
	\alpha^{(\ell)}_2={-2},\quad
	\alpha^{(\ell)}_3= 1, \\
\beta^{(\ell)}_1&=1,\quad
\beta^{(\ell)}_2=2,\quad
	\beta^{(\ell)}_3=1 ,\\
\alpha^{(\ell)}_{3k+4} &= 1,\quad
	\beta^{(\ell)}_{3k+4} = \frac{\beta^{(\ell-1)}_{k+2}}{\beta^{(\ell)}_{3k+3}\beta^{(\ell)}_{3k+2}},\quad
	\beta^{(\ell)}_{3k+5}=2-\beta^{(\ell)}_{3k+4},\\
	\alpha^{(\ell)}_{3k+5} &= -1- \frac{\alpha^{(\ell-1)}_{k+2}+1-\alpha^{(\ell)}_{3k+2}\beta^{(\ell)}_{3k+4}}{\beta^{(\ell)}_{3k+5}},\quad
	\alpha^{(\ell)}_{3k+6}=-1-\alpha^{(\ell)}_{3k+5},\\
	\beta^{(\ell)}_{3k+6} &= -1-\alpha^{(\ell)}_{3k+5} \alpha^{(\ell)}_{3k+6}.
\end{align*}

\begin{proposition} \label{TM3}
For $\ell \ge 0$ and $j=1, \ldots , 3^{\ell +1}$, we have

$\bullet$ $\nu_2(\alpha^{(\ell)}_{j})=-1$ if and only if
	$j=(3^{\ell+1}+1)/2$ or $j=(3^{\ell+1}+3)/2$; 
$\alpha^{(\ell)}_{j}\not= 0$ and $\nu_2(\alpha^{(\ell)}_{j})\geq 0$  otherwise;

$\bullet$ 
	$\nu_2(\beta^{(\ell)}_{2})=1$;
	$\nu_2(\beta^{(\ell)}_{j})=-2$ if and only if
	$j=(3^{\ell+1}+3)/2$;
	$\nu_2(\beta^{(\ell)}_{2})=0$ otherwise.

	As a conseqence, $\beta^{(\ell)}_j$ is nonzero for all $\ell\geq 0$ and $j=1,2,\ldots, 3^{\ell+1}$.
\end{proposition}

It follows from Proposition \ref{TM3} that the analogue of Theorem \ref{th:main} holds 
for the products $\tig_\ell(z)$, namely, there exists a positive real number $K$ such that, 
for every integer $b \ge 2$ and every integer $\ell \ge 2$, the inequality
$$
\Bigl| \prod_{h = 0}^\ell \, (1 - b^{-3^h} - b^{-2\cdot 3^h}) - {p \over q} \Bigr| 
> {1 \over q^2 \exp( K \log b \, \sqrt {\log q \, \log \log q})},
$$
holds for every rational number $p/q$ different from $\tig_\ell (b)$.

\begin{proof}
	We proceed by induction on $\ell,j$. We only display the more complicated steps. 
	We start with the $\alpha^{(\ell)}_{j}$'s. 
	Since $(3^{\ell+1}+1)/2$ is equal to $3k+5$ with $k=(3^\ell -3)/2$,
	\begin{align*}
		\nu_2(\alpha^{(\ell)}_{(3^{\ell+1}+1)/2})
		&=\nu_2\left(
	 -1- \frac{\alpha^{(\ell-1)}_{k+2}+1-\alpha^{(\ell)}_{3k+2}\beta^{(\ell)}_{3k+4}}{\beta^{(\ell)}_{3k+5}}
		\right)
=-1,\\
\end{align*}
because
$$
	\nu_2(\alpha^{(\ell-1)}_{k+2}) = -1,
	\quad
	\nu_2(\alpha^{(\ell)}_{3k+2})\geq 0,
	\quad
	\nu_2(\beta^{(\ell)}_{3k+4})=0,
	\quad
	\nu_2(\beta^{(\ell)}_{3k+5})=0.
$$
Then,
	\begin{align*}
		\nu_2(\alpha^{(\ell)}_{(3^{\ell+1}+3)/2})
		&=\nu_2\left(
		-1- \alpha^{(\ell)}_{(3^{\ell+1}+1)/2}\right)
=-1.
\end{align*}
	Since $(3^{\ell+1}+7)/2$ is equal to $3k+5$ with $k=(3^\ell -1)/2$,
	\begin{align*}
		\nu_2(\alpha^{(\ell)}_{(3^{\ell+1}+7)/2})
		&=\nu_2\left(
	 -1- \frac{\alpha^{(\ell-1)}_{k+2}+1-\alpha^{(\ell)}_{3k+2}\beta^{(\ell)}_{3k+4}}{\beta^{(\ell)}_{3k+5}}
		\right)
\geq 0,\\
\end{align*}
because
$$
	\nu_2(\alpha^{(\ell-1)}_{k+2}) = -1,
	\quad
	\nu_2(\alpha^{(\ell)}_{3k+2}) = -1,
	\quad
	\nu_2(\beta^{(\ell)}_{3k+4})=0,
	\quad
	\nu_2(\beta^{(\ell)}_{3k+5})=0.
$$
Then, 	\begin{align*}
		\nu_2(\alpha^{(\ell)}_{(3^{\ell+1}+9)/2})
		&=\nu_2\left(
		-1- \alpha^{(\ell)}_{(3^{\ell+1}+7)/2}\right)
\geq 0.
\end{align*}

Suppose that $j\not=  (3^{\ell+1}+1)/2, (3^{\ell+1}+3)/2, (3^{\ell+1}+7)/2, (3^{\ell+1}+9)/2$.
	If $j=3k+4$, then $\nu_2(\alpha^{(\ell)}_j)=0$.
	If $j=3k+5$ with $k\not=(3^\ell-3)/2, (3^\ell-1)/2$, then 
	\begin{align*}
		\nu_2(\alpha^{(\ell)}_{3k+5})
		&=\nu_2\left(
	 -1- \frac{\alpha^{(\ell-1)}_{k+2}+1-\alpha^{(\ell)}_{3k+2}\beta^{(\ell)}_{3k+4}}{\beta^{(\ell)}_{3k+5}}
		\right)
\geq 0,\\
\end{align*}
because
$$
	\nu_2(\alpha^{(\ell-1)}_{k+2}) \geq 0,
	\quad
	\nu_2(\alpha^{(\ell)}_{3k+2})\geq 0,
	\quad
	\nu_2(\beta^{(\ell)}_{3k+4})=0,
	\quad
	\nu_2(\beta^{(\ell)}_{3k+5})=0.
$$
Then, 
$$
\nu_2(\alpha^{(\ell)}_{3k+6})
=\nu_2(
	-1-\alpha^{(\ell)}_{3k+5}
)
\geq 0.
$$

We now deal with the $\beta^{(\ell)}_{j}$'s.
Since $(3^{\ell+1}+3)/2=3k+6$ with $k=(3^\ell-3)/2$,
$$
\beta^{(\ell)}_{(3^{\ell+1}+3)/2}
	=\nu_2( -1-\alpha^{(\ell)}_{3k+5} \alpha^{(\ell)}_{3k+6}) = -2,
$$
because
$$
	\nu_2( \alpha^{(\ell)}_{3k+5}) = \nu_2( \alpha^{(\ell)}_{3k+6}) = -1.
$$
Suppose that $j\not = (3^{\ell+1}+3)/2$. If $j=3k+4$ with $k\not=(3^{\ell+1}-1)/2$, then 
$$
\nu_2(\beta^{(\ell)}_{3k+4}) =\nu_2\left( \frac{\beta^{(\ell-1)}_{k+2}}{\beta^{(\ell)}_{3k+3}\beta^{(\ell)}_{3k+2}}\right) =0,
$$
because
$$
\nu_2(\beta^{(\ell-1)}_{k+2})=0,  
\nu_2(\beta^{(\ell)}_{3k+3})=0,
\nu_2(\beta^{(\ell)}_{3k+2})=0.
$$
If $j=3k+4$ with $k=(3^{\ell+1}-1)/2$, then 
$$
\nu_2(\beta^{(\ell)}_{3k+4}) =\nu_2\left( \frac{\beta^{(\ell-1)}_{k+2}}{\beta^{(\ell)}_{3k+3}\beta^{(\ell)}_{3k+2}}\right) =0,
$$
because
$$
\nu_2(\beta^{(\ell-1)}_{k+2})=-2,  
\nu_2(\beta^{(\ell)}_{3k+3})=-2,
\nu_2(\beta^{(\ell)}_{3k+2})=0.
$$
In both case,
$$
\nu_2(	\beta^{(\ell)}_{3k+5})=\nu_2(2-\beta^{(\ell)}_{3k+4})=0.
$$

If $j=3k+6$ with $k\not=(3^\ell-3)/2$,
$$
\beta^{(\ell)}_{3k+6}
	=\nu_2( -1-\alpha^{(\ell)}_{3k+5} \alpha^{(\ell)}_{3k+6}) = 0,
$$
because
$$
	\nu_2( \alpha^{(\ell)}_{3k+5}) \geq 0 , \nu_2( \alpha^{(\ell)}_{3k+6}) \geq 0,
	$$
	and if $\nu_2( \alpha^{(\ell)}_{3k+5} )=0$, then $ \nu(\alpha^{(\ell)}_{3k+6}) \geq 1 $.

\end{proof}

\end{document}